\input amstex
\input amsppt.sty
\magnification=\magstep1
\hsize=30truecc
\vsize=22.2truecm
\baselineskip=16truept
\TagsOnRight
\pageno=1
\nologo
\def\Z{\Bbb Z}
\def\N{\Bbb N}

\def\l{\left}
\def\r{\right}
\def\bg{\bigg}
\def\({\bg(}
\def\[{\bg\lfloor}
\def\){\bg)}
\def\]{\bg\rfloor}
\def\t{\text}
\def\f{\frac}

\def\p{\ (\roman{mod}\ p)}

\def\se {\subseteq}

\def\bi{\binom}
\def\eq{\equiv}

\def\ls{\leqslant}
\def\gs{\geqslant}
\def\mo{\roman{mod}}
\def\lcm{\roman{lcm}}

\def\al{\alpha}
\def\da{\delta}

\def\M#1#2{\thickfracwithdelims[]\thickness0{#1}{#2}_m}

\def\Mm#1#2#3{\thickfracwithdelims[]\thickness0{#1}{#2}_{#3}}

\def\mm#1#2#3{\thickfracwithdelims\{\}\thickness0{#1}{#2}_{#3}}

\def\Proof{\noindent{\it Proof}}

\def\Remark{\medskip\noindent{\it  Remark}}

\def\Ack{\medskip\noindent {\bf Acknowledgments}}
\hbox {J. Number Theory 126(2007), no.\,2, 287--296.}
\bigskip
\topmatter
\title Congruences for Sums of Binomial Coefficients\endtitle
\author Zhi-Wei Sun$^1$ and Roberto Tauraso$^2$\endauthor
\leftheadtext{Zhi-Wei Sun and Roberto Tauraso}
\affil $^1$Department of Mathematics, Nanjing University\\
 Nanjing 210093, People's Republic of China\\zwsun\@nju.edu.cn
\\{\tt http://math.nju.edu.cn/$^\sim$zwsun}
\medskip
$^2$Dipartimento di Matematica
\\Universit\`a di Roma ``Tor Vergata"
\\Roma 00133, Italy
\\tauraso\@mat.uniroma2.it
\\{\tt http://www.mat.uniroma2.it/$\sim$tauraso}
\endaffil
\abstract
Let $q>1$ and $m>0$ be relatively prime
integers. We find an explicit period $\nu_m(q)$
such that for any integers $n>0$ and $r$ we have
$$\Mm{n+\nu_m(q)}r{m}(a)\eq\Mm{n}r{m}(a)\ \ (\mo\ q)$$
whenever $a$ is an integer with $\gcd(1-(-a)^m,q)=1$,
or $a\eq-1\ (\mo\ q)$, or $a\eq1\ (\mo\ q)$ and $2\mid m$,
where $\Mm nrm(a)=\sum_{k\eq r\ (\mo\ m)}\bi nka^k$.
This is a further extension of a congruence of Glaisher.
\endabstract
\thanks 2000 {\it Mathematics Subject Classifications}.\,Primary 11B65;
Secondary 05A10,\,11A07.\newline\indent The first author is
supported by the National Science Fund for Distinguished Young
Scholars (Grant No. 10425103) in China.
\endthanks
\endtopmatter
\document

\heading{1. Introduction and main results}\endheading

Let $\N=\{0,1,2,\ldots\}$ and $\Z^+=\{1,2,3,\ldots\}$.
Following [S95, S02], for $m\in\Z^+$, $n\in\N$ and $r\in\Z$
we set
$$\M{n}r=\sum\Sb 0\ls k\ls n\\k\eq r\,(\mo\ m)\endSb \bi nk
=|\{X\se\{1,\ldots,n\}:\,|X|\eq r\ (\mo\ m)\}|\tag1.1$$
and
$$\mm nrm=\sum\Sb 0\ls k\ls n\\k\eq r\,(\mo\ m)\endSb(-1)^{\f{k-r}m}\bi nk
  =\Mm nr{2m}-\Mm n{r+m}{2m}.$$
Such sums occur in several topics of number theory or combinatorics.
(See, e.g., [SS, H, GS, S02].)

Let $p$ be an odd prime.
In 1899 J. W. L.
Glaisher obtained the
following congruence:
$$\Mm{n+p-1}r{p-1}\eq\Mm nr{p-1}\ (\mo\ p)
\ \ \t{for any}\ n\in\Z^+\ \t{and}\ r\in\Z.$$
Since an odd integer is not divisible by $p-1$,
this implies Hermite's result that $\Mm n0{p-1}\eq 1\ (\mo\ p)$
for $n=1,3,5,\ldots$ (cf. L. E. Dickson [D, p.\,271]).
A sophisticated proof of Glaisher's congruence can be found in
A. Granville [G97]; the first author observed in 2004 that Glaisher's congruence
can be proved immediately by induction on $n$.

   Before stating our further extension of
Glaisher's result, let us introduce some notations.

  Let $m\in\Z^+$, $n\in\N$ and $r\in\Z$. We set
  $$\Mm nrm(a)=\sum\Sb 0\ls k\ls n\\k\eq r\,(\mo\ m)\endSb\bi nka^k
  \quad \ \t{for}\ a\in\Z.\tag1.2$$
  Obviously $\Mm nrm(1)=\Mm nrm$, and
  $$\Mm nrm(-1)=\cases(-1)^r\Mm nrm&\t{if}\ 2\mid m,
  \\(-1)^r\mm nrm&\t{if}\ 2\nmid m.\endcases$$
  It is easy to see that
  $$\Mm{n+1}rm(a)=\Mm nrm(a)+a\Mm n{r-1}m(a).\tag1.3$$

   Let $a,b\in\Z$ and $q,m,n\in\Z^+$. Clearly
   $$\align&(x+a)^n\eq x^n+b\quad \mo\ (q,x^m-1)
   \\\iff&\sum_{r=0}^{m-1}\sum\Sb 0\ls k<n\\k\eq r\,(\mo\ m)\endSb
   \bi nk x^ka^{n-k}\eq b
   \quad \mo \ (q,x^m-1)
   \\\iff&\sum\Sb 0\ls k<n\\k\eq r\,(\mo\ m)\endSb\bi nk a^{n-k}
   \eq \cases b\ (\mo\ q)&\t{if}\ r=0,\\0\ (\mo\ q)&\t{if}\ 0<r<m\endcases
   \\\iff&\sum\Sb 1\ls k\ls n\\k\eq r\,(\mo\ m)\endSb\bi nka^k
   \eq\cases b\ (\mo\ q)&\t{if}\ r\eq n\ (\mo\ m),\\0\ (\mo\ q)&\t{otherwise}.\endcases
   \endalign$$
   (See also [G05].)
   Now that the congruence condition $(x+a)^n\eq x^n+a\ \mo\ (n,x^m-1)$
   plays a central role in the polynomial time primality test
   given by Agrawal, Kayal and Saxena [AKS], it is interesting to
   investigate periodicity of $\Mm nrm(a)$ mod $q$ (where $r\in\Z$)
   with respect to $n$.

  Let $q>1$ and $m>0$ be integers with $\gcd(q,m)=1$,
  where $\gcd(q,m)$ denotes the greatest common divisor of $q$ and $m$. Write $q$ in the
  factorization form $\prod_{s=1}^tp_s^{\al_s}$ where $p_1,\ldots,p_t$
  are distinct primes and $\al_1,\ldots,\al_t\in\Z^+$. We define
  $$\nu_m(q)=\lcm[p_1^{\al_1-1}(p_1^{\beta_1}-1),
  \ldots,p_t^{\al_t-1}(p_t^{\beta_t}-1)],\tag1.4$$
  where $\lcm[n_1,\ldots,n_t]$ represents the least common
  multiple of those $n_s\in\Z^+$ with $1\ls s\ls t$, and each $\beta_s$
  is the order of $p_s$ modulo $m$ (i.e., $\beta_s$
  is the smallest positive integer
  with $p_s^{\beta_s}\eq1\ (\mo\ m)$).
  Clearly $\nu_1(q)=\lcm[\varphi(p_1^{\al_1}),\ldots,\varphi(p_t^{\al_t})]$
  divides $\varphi(q)$, where $\varphi$ is Euler's totient function.
  Since $\varphi(p_s^{\al_s})\mid \nu_m(q)$ for each $s=1,\ldots,t$,
  if $a\in\Z$ is relatively prime to $q$, then by Euler's theorem
  $a^{\nu_m(q)}\eq1\ (\mo\ p_s^{\al_s})$ and therefore
  $a^{\nu_m(q)}\eq1\ (\mo\ q)$.
  Note also that $\nu_{p-1}(p^{\al})=\varphi(p^{\al})$
  for any prime $p$ and $\al\in\Z^+$.

  Now we present our first theorem.

\proclaim{Theorem 1.1} Let $q>1$ and $m>0$ be integers
with $\gcd(q,m)=1$. Let $T\in\Z^+$ be a multiple of $\nu_m(q)$,
and let $l\in\N$, $n\in\Z^+$ and $r\in\Z$.
Then
$$\sum_{k=0}^n(-1)^k\bi nk\Mm{kT+l}rm
\eq\cases 2^l(1-2^T)^n/m\ (\mo\ q^n)&\t{if}\ 2\nmid m,
\\\da_{l,0}(-1)^r/m\ (\mo\ q^n)&\t{if}\ 2\mid m,\endcases\tag1.5$$
where the Kronecker symbol $\da_{l,0}$ takes $1$ or $0$
according as $l=0$ or not.
\endproclaim

Actually Theorem 1.1 is implied by the following more general result
whose proof will be given in Section 2.

\proclaim{Theorem 1.2} Let $q>1$ be an integer relatively prime to both $m\in\Z^+$
and $\sum_{j=0}^{m-1}(-a)^j$ where $a\in\Z$. Let $l\in\N$ and $r\in\Z$.
If $n,T\in\Z^+$ and $\nu_m(q)\mid T$, then we have
$$\sum_{k=0}^n(-1)^k\bi nk\Mm{kT+l}r{m}(a)\eq\f{(a+1)^l}m\l(1-(a+1)^T\r)^n\ (\mo\ q^n).\tag1.6$$
\endproclaim

Now we explain why Theorem 1.1 follows from Theorem 1.2.
In the case $2\nmid m$, since $\sum_{j=0}^{m-1}(-1)^j=1$ we have (1.5)
by applying Theorem 1.2 with $a=1$.
In the case $2\mid m$, (1.5) also holds because
$$(-1)^r\Mm{kT+l}r{m}=(-1)^r\Mm{kT+l}r{m}(1)=\Mm{kT+l}r{m}(-1)$$
and therefore
$$\align&(-1)^r\sum_{k=0}^n(-1)^{k}\bi nk\Mm{kT+l}rm
\\=&\sum_{k=0}^n(-1)^k\bi nk\Mm{kT+l}rm(-1)
\eq\f{\da_{l,0}}m\ (\mo\ q^n)\endalign$$
with the help of Theorem 1.2 in the case $a=-1$.

 \proclaim{Corollary 1.3} Let $q>1$ and $m>0$ be integers with $\gcd(q,m)=1$.
 And let $l\in\N$, $n\in\Z^+$ and $r\in\Z$.

 {\rm (i)} Let $a$ be any integer with $\gcd(q,\sum_{j=0}^{m-1}(-a)^j)=1$.  Then
 $$\Mm{l+\nu_m(q)}rm(a)-\Mm lrm(a)
 \eq\cases0\ (\mo\ q_0),\\-(a+1)^l/m\ (\mo\ q/q_0),\endcases
 \tag1.7$$
 where $q_0$ is the largest divisor of $q$ relatively prime to $a+1$.
 Moreover, for each $k=1,2,3,\ldots$ we have
 $$\aligned&\Mm{k\nu_m(q)+l}r{m}(a)-\sum_{j=0}^{n-1}(-1)^{n-1-j}\bi{k-1-j}{n-1-j}\bi kj
 \Mm{j\nu_m(q)+l}r{m}(a)
 \\&\qquad\eq\f{(a+1)^l}{m}\sum_{n\ls j\ls k}\bi kj\l((a+1)^{\nu_m(q)}-1\r)^j\ (\mo\ q^{n}).
 \endaligned\tag1.8$$

 {\rm (ii)} Suppose that $m$ is even. For any $k\in\Z^+$ we have
 $$\aligned&\Mm{k\nu_m(q)+l}r{m}-\sum_{j=0}^{n-1}(-1)^{n-1-j}\bi{k-1-j}{n-1-j}\bi kj
 \Mm{j\nu_m(q)+l}r{m}
 \\&\qquad\quad\eq\da_{l,0}\f{(-1)^{n+r}}m\bi{k-1}{n-1}\ (\mo\ q^n).
 \endaligned\tag1.9$$
 In particular,
  $$\Mm{l+\nu_m(q)}rm-\Mm lrm\eq\da_{l,0}\f{(-1)^{r-1}}m\ (\mo\ q).\tag1.10$$
  \endproclaim
 \Proof. (i) Suppose that $p^{\al}\|q$ (i.e., $p^{\al}\mid q$ but
 $p^{\al+1}\nmid q$) where $p$ is a prime and $\al\in\Z^+$.
 If $p\mid a+1$, then $p^{\al}|(a+1)^{\nu_m(q)}$
 since $\nu_m(q)\gs p^{\al-1}\gs\al$;
 if $p\nmid a+1$, then $(a+1)^{\nu_m(q)}\eq1\ (\mo\ p^{\al})$
 as $\varphi(p^\al)\mid \nu_m(q)$.
 Therefore (1.7) follows from (1.6) in the case $n=1$ and $T=\nu_m(q)$.
 Note that $(a+1)^l\eq0\ (\mo\ q/q_0)$ if $l$ is sufficiently large.

  Let $k\in\Z^+$. By Lemma 2.1 of [Su],
 $$a_k-\sum_{j=0}^{n-1}(-1)^{n-1-j}\bi{k-1-j}{n-1-j}\bi kja_j
 =\sum_{n\ls j\ls k}\bi kj(-1)^j\sum_{i=0}^j\bi ji(-1)^ia_i$$
 for any sequence $a_0,a_1,\ldots$ of complex numbers. Applying this we
 immediately obtain (1.8) by noting that
$$\sum_{i=0}^j\bi ji(-1)^i\Mm {i\nu_m(q)+l}rm(a)
\eq\f{(a+1)^l}m\l(1-(a+1)^{\nu_m(q)}\r)^j\ (\mo\ q^j)$$
in view of (1.6).

 (ii) Applying (1.8) with $a=-1$, we find that
 $$\Mm{k\nu_m(q)+l}r{m}(-1)-\sum_{j=0}^{n-1}(-1)^{n-1-j}\bi{k-1-j}{n-1-j}\bi kj
 \Mm{j\nu_m(q)+l}r{m}(-1)$$
 is congruent to $\da_{l,0}m^{-1}\sum_{n\ls j\ls k}\bi kj(-1)^j$ modulo $q^n$.
 Observe that
 $$\align\sum_{n\ls j\ls k}\bi kj(-1)^j=&\sum_{n\ls j\ls k}\(\bi{k-1}j(-1)^j-\bi{k-1}{j-1}(-1)^{j-1}\)
 \\=&\bi{k-1}k(-1)^k-\bi{k-1}{n-1}(-1)^{n-1}=(-1)^n\bi{k-1}{n-1}.\endalign$$
 As $2\mid m$, we also have
 $$\Mm{j\nu_m(q)+l}rm(-1)=(-1)^r\Mm{j\nu_m(q)+l}rm\quad\t{for}\ j=0,1,2,\ldots.$$
 So (1.9) follows.  In the case $k=n=1$, (1.9) yields (1.10). We are done. \qed

 \Remark\ 1.1. Let $q>1$ and $m>0$ be relatively prime integers.
 Let $a$ be an integer such that $\gcd(1-(-a)^m,q)=1$, or $a\eq-1\ (\mo\ q)$, or $a\eq1\ (\mo\ q)$ and $2\mid m$.
 By Corollary 1.3(i), we have the following extension of Glaisher's periodic result:
$$\Mm{n+\nu_m(q)}rm(a)\eq\Mm nrm(a)\ \ (\mo\ q)\ \ \t{for any}\ n\in\Z^+\ \t{and}\ r\in\Z.\tag1.11$$
(Note that $\Mm nrm(-a)=(-1)^r\Mm nrm(a)$ if $2\mid m$.)

 \proclaim{Corollary 1.4}  Let $q>1$ be an integer relatively prime to $m\in\Z^+$.
 And let $k\in\Z^+$, $l\in\N$ and $r\in\Z$. Then
$$\aligned &\Mm{k\nu_m(q)+l}rm-k\Mm{\nu_m(q)+l}rm+(k-1)\Mm lrm
\\\eq&\cases\da_{l,0}(-1)^r(k-1)/m\ (\mo\ q^2)&\t{if}\ 2\mid m,
\\2^l(2^{k\nu_m(q)}-1-k(2^{\nu_m(q)}-1))/m\ (\mo\ q^2)&\t{if}\ 2\nmid m.
\endcases\endaligned\tag1.12$$
\endproclaim
\Proof.  In the case $2\mid m$, we get the desired congruence by applying (1.9) with $n=2$.
When $2\nmid m$, putting $a=1$ in (1.8) we obtain
$$\align&\Mm{k\nu_m(q)+l}rm-k\Mm{\nu_m(q)+l}rm+(k-1)\Mm lrm
\\\eq&\f{2^l}m\sum_{2\ls j\ls k}\bi kj(2^{\nu_m(q)}-1)^j
=\f{2^l}m(2^{k\nu_m(q)}-1-k(2^{\nu_m(q)}-1))\ (\mo\ q^2).
\endalign$$
This completes the proof. \qed

\Remark\ 1.2. Let $p$ be an odd prime. Let $k\in\Z^+$ and $r\in\{0,1,\ldots,p-2\}$.
As $\nu_{p-1}(p)=p-1$, by Corollary 1.4 we have
$$\Mm{k(p-1)}r{p-1}\eq k\Mm{p-1}r{p-1}-(k-1)\Mm 0r{p-1}+(-1)^r\f{k-1}{p-1}\ (\mo\ p^2).$$
As $0\ls r<p-1$ and $1/(p-1)\eq-p-1\ (\mo\ p^2)$, this turns out to be
 $$\Mm{k(p-1)}r{p-1}\eq k\bi{p-1}r-(-1)^r(k-1)(p+1)+\da_{r,0}\
 (\mo\ p^2).\tag1.13$$
In the case $r=0$, this solves a problem
 proposed by V. Dimitrov [Di].
 \medskip

    Let $p$ be any odd prime and let $\al,n\in\Z^+$.
 As $\nu_{p-1}(p^\al)=p^\al-p^{\al-1}$, by Remark 1.1 we have
  $$\Mm{p^{\al}n}r{p-1}\eq
  \Mm{p^{\al-1}n}r{p-1}\ (\mo\ p^{\al})\ \ \t{for any}\ r\in\Z.\tag1.14$$
  In 1953, by using some deep properties of Bernoulli numbers,
 L. Carlitz [C] extended Hermite's congruence in the
  following way:
  $$p+(p-1)\sum\Sb 0<k<p^{\al-1}n\\p-1\mid k\endSb\bi{p^{\al-1}n}k\eq0
  \ \ (\mo\ p^{\al}).$$
  When $p-1\mid n$, this follows from (1.10), for,
  $\nu_{p-1}(p^{\al})$ divides $p^{\al-1}n$ and hence
  $$\Mm{p^{\al-1}n}0{p-1}\eq\Mm 00{p-1}-\f1{p-1}=1-\f1{p-1}
  \ (\mo\ p^{\al}).$$

 Let $q>1$ and $m>0$ be integers with $\gcd(q,m)=1$.
 Let $a$ be an integer with $\gcd(1-(-a)^m,q)=1$, or $a\eq-1\ (\mo\ q)$, or $a\eq1\ (\mo\ q)$ and $2\mid m$.
 What is the smallest positive integer $\mu_m(a,q)$ such that
 $$\Mm{n+\mu_m(a,q)}rm(a)\eq\Mm nrm(a)\ \ (\mo\ q)\tag1.15$$
 holds for all $n\in\Z^+$ and $r\in\Z$?
 Clearly $\mu_m(0,q)=1$, and
 $\mu_m(a,q)|\nu_m(q)$ by (1.11).
 (If $\mu_m(a,q)\nmid\nu_m(q)$, then the least positive residue
 of $\nu_m(q)$ mod $\mu_m(a,q)$ would be
 a period smaller than $\mu_m(a,q)$.)

 \proclaim{Conjecture 1.5} Let $q>1$ and $m>0$ be
 integers with $\gcd(q,m)=1$ and $q\not\eq0\ (\mo\ 3)$.
 Then $\nu_m(q)$ is the maximal value of
 $\mu_m(a,q)$, where $a$ is an integer with $\gcd(1-(-a)^m,q)=1$, or $a\eq-1\ (\mo\ q)$,
 or $a\eq 1\ (\mo\ q)$ and $2\mid m$.
 \endproclaim

  Now we give an example to illustrate our conjecture.
  \medskip

 {\it Example} 1.1. (i) Since the order of 3 modulo 7 is 6, we have
 $\nu_7(9)=3(3^6-1)=2184$. For any given $a\in\Z$, clearly
 $$1-(-a)^7=1+a^3a^3a\eq1+a^3\eq1+a\ (\mo\ 3)$$
 since $a^3\eq a\ (\mo\ 3)$, thus $\gcd(1-(-a)^7,9)=1$
 if and only if $a\not\eq2\ (\mo\ 3)$.
 Through computation we obtain that
 $$\mu_7(-1,9)=1092,\ \mu_7(1,9)=\mu_7(-2,9)=\mu_7(4,9)=546,
  \ \mu_7(\pm3,9)=3.$$

  (ii) The order of 5 modulo 7 is 6, thus
 $\nu_7(5)=5^6-1=15624$. For any given $a\in\Z$, clearly
 $1-(-a)^7=1+a^5a^2\eq1+a^3\ (\mo\ 5)$,
  thus $5\nmid 1-(-a)^7$
 if and only if $a\not\eq-1\ (\mo\ 5)$.
 By computation we find that
 $$\mu_7(1,5)=868,\ \mu_7(-1,5)=1736,\ \mu_7(2,5)=2232,
 \ \mu_7(-2,5)=15624.$$

 (iii) Clearly $\nu_6(11)=11^2-1=120$. By computation, $\mu_6(\pm1,11)=60$
 and $\mu_6(a,11)=120$ for any integer $a\not\eq0,\pm1\ (\mo\ 11)$.
 Note that $4(a^4+a^2+1)=(2a^2+1)^2+3\not\eq0\ (\mo\ 11)$
 since $-3$ is a quadratic non-residue modulo $11$.
 Thus, if $a\not\eq\pm1\ (\mo\ 11)$ then
 $1-(-a)^6=(1-a^2)(a^4+a^2+1)$ is relatively prime to $11$.

 \heading{2. Proof of Theorem 1.2}\endheading

 In this section we work with congruences in the ring of algebraic integers.
 The reader may consult [IR, pp.\,66--69]
 for the basic knowledge of algebraic integers.

 \proclaim{Lemma 2.1}
 Let $a\in\Z$ and $m\in\Z^+$, and let
 $q>1$ be an integer relatively prime to $m\sum_{j=0}^{m-1}(-a)^j$.
 If $\zeta\not=1$ is an $m$-th root of unity, then we have the congruence
 $$(1+a\zeta)^{\nu_m(q)}\eq1\ \ (\mo\ q)\tag2.1$$
 in the ring of algebraic integers.
 \endproclaim
 \Proof. Let $p$ be any prime divisor of $q$, and let $\beta$ be the order
 of $p$ modulo $m$.
  Below we use induction to show that
  $$(1+a\zeta)^{p^{\al-1}(p^{\beta}-1)}\eq1\ \ (\mo\ p^{\al})\tag2.2$$
  for every $\al=1,2,3,\ldots$.

 Since $p\mid \bi pk$ for $k=1,\ldots,p-1$ and $a^p\eq a\ (\mo\ p)$
 by Fermat's little theorem, we have
 $$(1+a\zeta)^p=1+a^p\zeta^p+\sum_{k=1}^{p-1}\bi pk(a\zeta)^k
 \eq1+a\zeta^p\ (\mo\ p),$$
 hence
 $$(1+a\zeta)^{p^2}\eq(1+a\zeta^p)^p\eq1+a\zeta^{p^2}\ (\mo\ p)$$
 and so on. Thus
 $$(1+a\zeta)^{p^{\beta}}\eq1+a\zeta^{p^{\beta}}=1+a\zeta\ (\mo\ p).$$
 (Recall that $p^{\beta}\eq1\ (\mo\ m)$ and $\zeta^m=1$.)
 Clearly
 $$\align\prod_{0<j<m}\f{1+ae^{2\pi ij/m}}{-e^{2\pi ij/m}}
 =&\prod_{0<j<m}(x-e^{-2\pi ij/m})\bg|_{x=-a}
 \\=&\lim_{x\to -a}\f{x^m-1}{x-1}=\sum_{j=0}^{m-1}(-a)^j
 \endalign$$
 and so $z=1+a\zeta$ divides $c=\sum_{j=0}^{m-1}(-a)^j$
 in the ring of algebraic
 integers.
 Therefore
 $$cz^{p^{\beta}-1}\eq\f{c}z z^{p^{\beta}}\eq\f{c}zz\eq c\ (\mo\ p)$$
 and hence $z^{p^{\beta}-1}\eq1\ (\mo\ p)$ since $p\nmid c$.
 This proves (2.2) in the case $\al=1$.

   Now let $\al\in\Z^+$ and suppose that (2.2) holds. Then
   $z^{p^{\al-1}(p^{\beta}-1)}=1+p^{\al}\omega$ for some
   algebraic integer $\omega$. It follows that
  $$z^{p^{\al}(p^{\beta}-1)}=(1+p^{\al}\omega)^p
  \eq1+\bi p1 p^{\al}\omega\eq1\ (\mo\ p^{\al+1}).$$
  This concludes the induction step.

   For any $q_1,q_2\in\Z$ with $\gcd(q_1,q_2)=1$, there
   are $x_1,x_2\in\Z$ such that $q_1x_1+q_2x_2=1$,
   If an algebraic integer $\omega$ is divisible
   by both $q_1$ and $q_2$, then $\omega=q_1(\omega x_1)+q_2(\omega x_2)$
   is divisible by $q_1q_2$ in the ring of algebraic integers.
   Therefore (2.1) is valid in view of what we have proved. \qed

   \Remark\ 2.1. Write an integer $q>1$ in the form
   $p_1^{\al_1}\cdots p_t^{\al_t}$,
   where $p_1,\ldots,p_t$ are distinct primes and $\al_1,\ldots,\al_t\in\Z^+$.
   Let $m$ be a positive integer dividing
   $p_s-1$ for all $s=1,\ldots,t$.
   And let $g$ be an integer with $g\eq g_s^{\varphi(p_s^{\al_s})/m}
   \ (\mo\ p_s^{\al_s})$
   for $s=1,\ldots,t$, where
   $g_s$ is a primitive root modulo $p_s$.
    Clearly $g^m\eq1\ (\mo\ q)$.
   Suppose that $j\in\Z^+$ and $j<m$. Then $p_s-1\nmid j\varphi(p_s^{\al_s})/m$
   and hence $g^j\not\eq1\ (\mo\ p_s)$. Therefore $\gcd(g^j-1,q)=1$.
   If $p_s\mid 1+ag^j$, then $-a\eq g^{m-j}\not\eq1\ (\mo\ p_s)$
   but $(a+1)\sum_{i=0}^{m-1}(-a)^i=1-(-a)^m\eq1-g^{(m-j)m}\eq0\ (\mo\ p_s).$
   Thus, if $\gcd(\sum_{i=0}^{m-1}(-a)^i,q)=1$, then
   $\gcd(1+ag^j,q)=1$, and hence
   $$(1+ag^j)^{\nu_m(q)}\eq1\ (\mo\ q)\tag2.3$$
   which is an analogue of (2.1).

  \noindent{\it Proof of Theorem 1.2}.
   Set $\zeta=e^{2\pi i/m}$.
   For any $h\in\Z$, we clearly have
   $$\sum_{j=0}^{m-1}\zeta^{jh}=\cases m&\t{if}\ m\mid h,
   \\0&\t{otherwise}.\endcases$$
  If $n\in\N$ then
  $$\align m\Mm nrm(a)
  =&\sum_{k=0}^n\bi nka^k\sum_{j=0}^{m-1}\zeta^{j(k-r)}
  \\=&\sum_{j=0}^{m-1}\zeta^{-jr}\sum_{k=0}^n\bi nka^k\zeta^{jk}
  =\sum_{j=0}^{m-1}\zeta^{-jr}(1+a\zeta^j)^n.
  \endalign$$

  Now let $T\in\Z^+$ be a multiple of $\nu_m(q)$, and fix a positive
  integer $n$.
  By the above,
  $$\align&m\sum_{k=0}^n(-1)^k\bi nk\Mm{kT+l}rm(a)
  \\=&\sum_{j=0}^{m-1}\zeta^{-jr}
  \sum_{k=0}^n\bi nk(-1)^k(1+a\zeta^j)^{kT+l}
  \\=&\sum_{j=0}^{m-1}\zeta^{-jr}(1+a\zeta^j)^l\l(1-(1+a\zeta^j)^{T}\r)^n
  \\\eq&(1+a)^l\l(1-(1+a)^T\r)^n\ (\mo\ q^{n})
  \endalign$$
  where we have applied Lemma 2.1.
  This concludes our proof. \qed

  \Remark\ 2.2. Let $a,r\in\Z$ and $m\in\Z^+$, and let $q>1$
  be an integer relatively prime to $\sum_{j=0}^{m-1}(-a)^j$.
  Suppose that $m\mid p-1$ for any prime divisor $p$ of $q$.
  Obviously $\gcd(m,q)=1$. Choose $g\in\Z$ as in Remark 2.1. Then
   $g^m\eq1\ (\mo\ q)$, and for each $0<j<m$ we have
   $\gcd(g^j-1,q)=1$ as well as $(2.3)$.
   By modifying the proof of Theorem 1.2 slightly,
   we find that
  $$m\Mm nrm(a)\eq\sum_{j=0}^{m-1}g^{-jr}(1+ag^j)^n
  =(a+1)^n+\sum_{0<j<m}g^{-jr}a_j^n\ \ (\mo\ q)$$
  for every $n\in\N$,
  where $a_j=1+ag^j\ (0<j<m)$ are relatively prime to $q$.
  If $q\mid a+1$ or $\gcd(a+1,q)=1$,
  then the function $f:\Z^+\to\Z$ given by $f(n)=\Mm nrm(a)$
  is $q$-normal in the sense that
  $$f(n)\eq\sum\Sb 1\ls j<q\\\gcd(j,q)=1\endSb c_jj^n\ (\mo\ q)
  \quad\ \t{for all}\ n\in\Z^+,\tag2.4$$
  where $c_j\ (1\ls j<q\ \t{and}\ \gcd(j,q)=1)$ are suitable integers.
  The concept of $q$-normal function was first introduced
  by Sun [S03] where the reader can find some $q$-normal functions
  involving Bernoulli polynomials.

 \Ack.  The joint work was done during the first author's stay
at Univ. Lyon-I as a visiting professor, thus he
is indebted to Prof. J. Zeng for the invitation and
hospitality. The authors also thank the referee for his/her helpful comments.

\enddocument